\documentclass[english,12pt]{article}
\usepackage{amsfonts}
\usepackage{amssymb}
\usepackage{amssymb, amsthm, amsmath}
\usepackage[latin1]{inputenc}
\usepackage{babel}
\usepackage{graphicx}
\usepackage{psfrag}
\usepackage{amstext}
\newtheorem{teorema}{\bf Theorem}

\newtheorem{observacao}{\bf Remark}

\newtheorem{definicao}{ \bf Definition}

\newcommand{\abs}[1]{\lvert #1 \rvert}

\begin{document}
\title{A Characterization of Minimal Surfaces in $S^{5}$ with Parallel Normal Vector Field}
\author{Rodrigo Ristow Montes
\thanks{ristow@mat.ufpb.br and ristow@math.wustl.edu}}
\date{Departamento de Matem\'atica , \\
      Universidade Federal da Para\'iba, \\[2mm]
      BR-- 58.051-900 ~ Jo\~ao Pessoa, P.B., Brazil \\}
\maketitle
\renewcommand{\thefootnote}{\fnsymbol{footnote}}
\renewcommand{\thefootnote}{\arabic{footnote}}
\setcounter{footnote}{0}
\thispagestyle{empty}
\begin{abstract}
\noindent
  In this paper we proof that the Holomorphic angle for compact minimal surfaces in the sphere $S^5$ with constant Contact angle and with a parallel normal vector field must be constant.
\end{abstract}

\smallskip
\noindent{\bf 2000 Math Subject Classification:} 53C42 - 53D10 - 53D35.
\section{Introduction}\mbox{}
 The notion of K\"ahler angle was introduced by  Chern and Wolfson
 in \cite{CW} and \cite{W}; it is a fundamental invariant
for  minimal surfaces in complex manifolds. Using the technique of
 moving frames, Wolfson obtained equations for the Laplacian and Gaussian
 curvature for an immersed minimal surface in $\mathbb{CP}^n$. Later, Kenmotsu
 in \cite{K}, Ohnita in \cite{O} and  Ogata in \cite{Og} classified  minimal surfaces
 with constant Gaussian curvature and constant  K\"ahler angle.\\
A few years ago, Li in  \cite{Z} gave a counterexample to the
conjecture of Bolton, Jensen and Rigoli (see \cite{BJRW}), according to which
a minimal  immersion (non-holomorphic, non anti-holomorphic, non totally real) of a  two-sphere in $\mathbb{CP}^n$
with constant  K\"ahler angle would have constant Gaussian curvature.\\
In  \cite{MV} we introduced the notion of contact angle, that can be
considered as  a new geometric invariant useful to investigate
the geometry of immersed surfaces in $S^{3}$. Geometrically, the contact angle
$(\beta)$ is the complementary angle between the contact
distribution and the tangent space of the surface. Also in
\cite{MV}, we deduced formulas for the Gaussian curvature and the Laplacian
of an immersed minimal surface in $S^3$, and we gave a
characterization of the  Clifford Torus as the only minimal
surface in $S^3$
with constant contact angle.\\
We define $\alpha$ to be the angle given by $ \cos \alpha = \langle
ie_1 , v \rangle \nonumber$, where $e_1$ and $v$ are defined in
section~\ref{sec:section2}. The holomorphic angle $\alpha$ is the analogue of
the K\"ahler
angle introduced by Chern and Wolfson in \cite{CW}.\\
Recently,  in \cite{RMV},  we construct a
family of minimal tori in $S^{5}$ with constant contact and
holomorphic angle. These tori are parametrized by the following circle
equation
\begin{equation}\label{eq:equacaoab}
a^2 + \left(b- \frac{\cos\beta}{1+\sin^2\beta}\right)^2 =  2\frac{\sin^4\beta}{(1+\sin^2\beta)^2},
\end{equation}
where  $a$ and $b$ are given in Section~\ref{parametro} (equation
(\ref{eq: segunda})).
In particular, when $a=0$ in \eqref{eq:equacaoab}, we recover the
examples found by Kenmotsu, in \cite{KK2}.  These examples are
defined for $0 < \beta < \frac{\pi}{2}$. Also, when $b=0$ in
\eqref{eq:equacaoab}, we find a new family of minimal tori in
$S^5$, and these tori are defined for $\frac{\pi}{4} < \beta <
\frac{\pi}{2}$. Also, in \cite{RMV}, when $\beta = \frac{\pi}{2}$, we
give an alternative proof of this classification of a Theorem from 
Blair in \cite{Blair}, and Yamaguchi, Kon and Miyahara in
\cite{YKM} for Legendrian minimal surfaces in $S^5$ with
constant Gaussian curvature.\\
In this paper, we classify minimal surfaces in $S^5$ with constant Contact angle and  with a  parallel normal vector field. We suppose that $e_3$ (in equation (\ref{eq:camposnormais})) is a parallel normal vector field, and we get the following
\begin{teorema}
The holomorphic angle  $( 0 < \alpha < \frac{\pi}{2})$ is constant for compact minimal surfaces in $S^5$ with constant Contact angle $\beta$ and null principal curvatures $a,b$ 
\end{teorema}


\section{Contact Angle for Immersed Surfaces in  $S^{2n+1}$}\label{sec:section2}
Consider in $\mathbb{C}^{n+1}$ the following objects:
\begin{itemize}
\item the Hermitian product: $(z,w)=\sum_{j=0}^n z^j\bar{w}^j$;
\item the inner product: $\langle z,w \rangle = Re (z,w)$;
\item the unit sphere: $S^{2n+1}=\big\{z\in\mathbb{C}^{n+1} | (z,z)=1\big\}$;
\item the \emph{Reeb} vector field in $S^{2n+1}$, given by: $\xi(z)=iz$;
\item the contact distribution in $S^{2n+1}$, which is orthogonal to $\xi$:
\[\Delta_z=\big\{v\in T_zS^{2n+1} | \langle \xi , v \rangle = 0\big\}.\]
\end{itemize}
We observe that $\Delta$ is invariant by the complex structure of $\mathbb{C}^{n+1}$.

Let now $S$ be an immersed orientable surface in $S^{2n+1}$.
\begin{definicao}
The \emph{contact angle} $\beta$ is the complementary angle between the
contact distribution $\Delta$ and the tangent space $TS$ of the
surface.
\end{definicao}
Let $(e_1,e_2)$ be  a local frame of $TS$, where $e_1\in
TS\cap\Delta$. Then $\cos \beta = \langle \xi , e_2 \rangle
$. Finally, let $v$ be the unit vector in the direction of  the orthogonal projection  of $e_2$ on $\Delta$,
defined by the following relation
\begin{eqnarray}\label{eq:campoe2}
e_2 = \sin\beta v + \cos\beta \xi.
\end{eqnarray}
\section{Equations for Gaussian curvature and Laplacian of a minimal surface in
  $S^5$}\label{parametro}
In this section, we deduce the equations for the Gaussian curvature  and for the Laplacian
of a minimal surface in $S^5$ in terms of the
contact angle and the holomorphic angle.
Consider the normal vector fields
\begin{eqnarray}\label{eq:camposnormais}
e_3                   & = & i \csc \alpha e_1 - \cot \alpha v \nonumber\\
e_4                   & = & \cot \alpha e_1 + i \csc \alpha v \\
e_5                   & = & \csc \beta \xi - \cot \beta e_2
\nonumber
\end{eqnarray}
where $\beta \neq {{0,\pi}}$ and  $ \alpha \neq{{0,\pi}}$.
We will call $(e_j)_{1 \leq \ j \leq 5}$ an  \emph{adapted frame}.

Using (\ref{eq:campoe2}) and (\ref{eq:camposnormais}), we get
\begin{eqnarray}\label{eq:camposdistribution}\
v  =  \sin \beta e_2  - \cos \beta e_5 , \quad  iv   =  \sin \alpha e_4 - \cos \alpha e_1 \\
\xi  =  \cos \beta e_2 + \sin \beta e_5 \nonumber
\end{eqnarray}
It follows from (\ref{eq:camposnormais}) and
(\ref{eq:camposdistribution}) that
\begin{eqnarray}\label{eq:inormais}
ie_1 & = & \cos \alpha \sin \beta e_2 + \sin\alpha e_3 -\cos\alpha \cos \beta
e_5  \\
ie_2 & = & - \cos \beta z - \cos\alpha \sin\beta e_1  + \sin\alpha \sin
\beta e_4 \nonumber
\end{eqnarray}
Consider now the dual basis $(\theta^j)$ of $(e_j)$.
The connection forms $(\theta_k^j)$ are given by
\begin{eqnarray}
De_j    =   \theta_j^k e_k, \nonumber
\end{eqnarray}
and the second fundamental form  with respect to this frame are given by
\begin{equation}
\begin{array}{lclll}
II^j      & = & \theta_1^j \theta^1  +  \theta_2^j \theta^2; \quad
j=3, ..., 5 \nonumber.
\end{array}
\end{equation}
Using (\ref{eq:inormais}) and differentiating $v$ and $\xi$ on the
surface $S$, we get
\begin{eqnarray}\label{eq:dif}
D\xi & = & -\cos\alpha \sin\beta \theta^2 e_1 + \cos \alpha \sin \beta
 \theta^1 e_2 + \sin \alpha \theta^1 e_3 + \sin\alpha \sin \beta \theta^2 e_4\nonumber\\
          && - \cos \alpha \cos \beta \theta^1 e_5,  \\
Dv   & = & (\sin \beta \theta_2^1 - \cos \beta \theta_5^1)e_1 +
 \cos\beta(d\beta-\theta_5^2)e_2 + ( \sin\beta \theta_2^3 - \cos \beta
 \theta_5^3)e_3 \nonumber \\
          && +  ( \sin\beta \theta_4^2 - \cos \beta \theta_5^4)e_4 + \sin
 \beta(d\beta + \theta_2^5)e_5 \nonumber.
\end{eqnarray}
Differentiating $e_3$, $e_4$ and  $e_5$, we have
\begin{eqnarray}\label{eq:intrin}
\theta_3^1 & = &  -\theta_1^3 \nonumber \\
\theta_3^2 & = &  \phantom{-}\sin \beta ( d\alpha + \theta_4^1) - \cos \beta \sin \alpha
\theta^1 \nonumber \\
\theta_3^4 & = & \phantom{-} \csc \beta \theta_1^2 - \cot \alpha ( \theta_1^3 + \csc \beta
\theta_2^4) \nonumber \\
\theta_3^5 & = & \phantom{-} \cot \beta \theta_2^3 - \csc \beta \sin \alpha \theta^1 \nonumber \\
\theta_4^1 & = & -d\alpha - \csc \beta \theta_2^3 + \sin \alpha \cot \beta \theta^1 \nonumber \\
\theta_4^2 & = & - \theta_2^4 \nonumber \\
\theta_4^3 & = & \phantom{-} \csc \beta \theta_2^1 + \cot \alpha ( \theta_1^3 + \csc \beta
\theta_2^4) \nonumber \\
\theta_4^5 & = & \phantom{-} \cot \beta \theta_2^4 - \sin \alpha \theta^2 \nonumber \\
\theta_5^1 & = & -\cos \alpha \theta^2 -  \cot \beta \theta_2^1 \nonumber \\
\theta_5^2 & = & \phantom{-} d\beta + \cos \alpha \theta^1 \\  \label{eq: intrin}
\theta_5^3 & = & -\cot \beta \theta_2^3 + \csc \beta \sin \alpha \theta^1 \nonumber\\
\theta_5^4 & = & -\cot \beta \theta_2^4 + \sin \alpha \theta^2 \nonumber
\end{eqnarray}
The conditions of minimality and of symmetry  are equivalent to the following
equations:
\begin{eqnarray}\label{eq:minimal}
\theta_1^\lambda \wedge \theta^1 +  \theta_2^\lambda \wedge
\theta^2  = 0 = \theta_1^\lambda \wedge \theta^2 -
\theta_2^\lambda \wedge \theta^1.
\end{eqnarray}
On the surface $S$, we consider
\begin{eqnarray}
\theta_1^3 & = & a\theta^1 + b\theta^2 \nonumber
\end{eqnarray}
It follows from (\ref{eq:minimal}) that
\begin{eqnarray}\label{eq: segunda}
\theta_1^3 & = & \phantom{-} a\theta^1 + b\theta^2\nonumber\\
\theta_2^3 & = & \phantom{-} b\theta^1 - a\theta^2\nonumber\\
\theta_1^4 & = & \phantom{-} d\alpha + (b \csc\beta - \sin\alpha
\cot\beta)\theta^1 - a \csc \beta \theta^2\nonumber \\
\theta_2^4 & = & \phantom{-} d\alpha \circ J - a \csc\beta \theta^1 - (b
\csc\beta - \sin\alpha \cot\beta)\theta^2 \\
\theta_1^5   &  =  & \phantom{-} d\beta \circ J  - \cos \alpha \theta^2\nonumber\\
\theta_2^5   &  =  &  - d\beta - \cos \alpha \theta^1\nonumber
\end{eqnarray}
where $J$ is the complex structure  of $S$ is given by $Je_1=e_2$ and $Je_2=-e_1$.
Moreover, the normal connection forms are given by:
\begin{eqnarray}\label{eq:normalconexaobeta}
\theta_3^4 & = & - \sec\beta d\beta \circ J - \cot\alpha \csc\beta d\alpha
\circ J + a \cot\alpha \cot^2\beta \theta^1 \nonumber \\
&&+ ( b  \cot\alpha \cot^2\beta -
\cos\alpha \cot\beta \csc\beta + 2 \sec\beta \cos \alpha) \theta^2 \nonumber \\
\theta_3^5 & = & \phantom{-} (b \cot\beta - \csc\beta \sin\alpha) \theta^1 - a
\cot\beta \theta^2 \\
\theta_4^5 & = & \phantom{-} \cot\beta (d\alpha \circ J) - a \cot\beta
\csc\beta \theta^1 + ( -b\csc\beta \cot\beta + \sin\alpha (\cot^2\beta -1))\theta^2,\nonumber
\end{eqnarray}
while the Gauss equation is equivalent to the equation:
\begin{equation}\label{eq: Gauss}
\begin{array}{lcl}
d\theta_2^1 + \theta_k^1 \wedge \theta_2^k & = & \theta^1 \wedge \theta^2.
\end{array}
\end{equation}
Therefore, using  equations  (\ref{eq: segunda}) and (\ref{eq:
Gauss}), we have
\begin{eqnarray}\label{eq:curvatura1}
K  & = &  1 - |\nabla \beta|^2 - 2 \cos \alpha \beta_1 - \cos^2 \alpha
-(1+ \csc^2 \beta )(a^2 + b^2)
\nonumber\\
 &&+ 2b\sin\alpha\csc\beta\cot\beta  + 2 \sin \alpha \cot \beta \alpha_1 - |\nabla \alpha|^2 \nonumber\\
   && + 2 a \csc \beta
 \alpha_2 - 2b  \csc \beta \alpha_1 - \sin^2 \alpha \cot^2 \beta \nonumber\\
&&  =    1 - (1+ csc^2 \beta)(a^2+b^2)-2b\csc\beta
(\alpha_1-\sin\alpha \cot\beta)+2a \csc\beta \alpha_2
  \\
&& - |\nabla \beta + \cos \alpha e_1|^2 - |\nabla \alpha - \sin \alpha \cot
  \beta e_1|^2
        \nonumber
\end{eqnarray}
Using (\ref{eq:intrin}) and  the complex structure of $S$, we get
\begin{equation} \label{eq:conex}
\begin{array}{lcl}
\theta_2^1  &  =  &  \tan\beta(d\beta\circ J-2\cos\alpha\theta^2)
\end{array}
\end{equation}
Differentiating (\ref{eq:conex}), we conclude that
\begin{eqnarray}
d\theta_2^1 & = & ( -(1 + \tan^2
\beta)\abs{\nabla\beta}^2-\tan\beta\Delta\beta
-2\cos\alpha(1+2\tan^{2}\beta)\beta_1\nonumber\\
    &&  +2\tan\beta\sin\alpha\alpha_1-4\tan^{2}\beta\cos^{2}\alpha) \theta^1
    \wedge \theta^2 \nonumber
\end{eqnarray}
where $\Delta = tr \nabla^2 $ is the Laplacian of $S$. The Gaussian curvature is therefore
given by:
\begin{eqnarray}\label{eq:curvatura2}
K & = &  -(1 + \tan^2 \beta)\abs{\nabla\beta}^2-\tan\beta\Delta\beta
-2\cos\alpha(1+2\tan^{2}\beta)\beta_1 \nonumber\\
    &&
    +2\tan\beta\sin\alpha\alpha_1-4\tan^{2}\beta\cos^{2}\alpha.
\end{eqnarray}
From  (\ref{eq:curvatura1}) and (\ref{eq:curvatura2}), we obtain
the following formula for the Laplacian of $S$:
\begin{eqnarray}\label{eq:lapla}
\tan\beta\Delta\beta &  = & (1+ \csc^2 \beta)(a^2+b^2)+2b\csc\beta
(\alpha_1-\sin\alpha \cot\beta)-2a \csc\beta \alpha_2 \nonumber \\
&& - \tan^2 \beta( |\nabla \beta + 2 \cos \alpha e_1|^2
- |\cot \beta \nabla\alpha + \sin \alpha ( 1- \cot^2 \beta)e_1|^2 )\nonumber \\
 && + \sin^2 \alpha ( 1 -\tan^2 \beta)
\end{eqnarray}
\section{Gauss-Codazzi-Ricci equations for a minimal surface in $S^5$ with constant Contact angle $\beta$}
In this section, we will compute Gauss-Codazzi-Ricci equations for a minimal surface in $S^{5}$ with constant Contact angle $\beta$.\\
Using the connection form (\ref{eq: segunda}) and 
(\ref{eq:normalconexaobeta}) in the Codazzi-Ricci equations, we have
\begin{eqnarray}
d\theta_1^3 + \theta_2^3 \wedge \theta_1^2 + \theta_4^3 \wedge \theta_1^4 + \theta_5^3 \wedge \theta_1^5  &  =  & 0 \nonumber
\end{eqnarray}
This implies that
\begin{eqnarray}\label{eq:Codazzi1}
&& (b_1-a_2) + (a^2+ b^2) \cot\alpha \csc\beta \cot^2 \beta - a \cot\alpha(\csc^2\beta + \cot^2\beta) \alpha_2 \\
&& + b(cot\alpha(\csc^2\beta + \cot^2\beta)\alpha_1 - \cos\alpha \cot\beta (\csc^2\beta + \cot^2\beta -3\sec^2\beta (1+\sin^2 \beta)))\nonumber\\
&& - \cos\alpha \csc \beta (2(\cot\beta - \tan\beta)\alpha_1  - \sin\alpha (\cot^2 \beta -3)) + \cot\alpha \csc\beta |\nabla \alpha|^2 =  0\nonumber
\end{eqnarray}
Replacing the following (\ref{eq:normalconexaobeta}) in the Codazzi-Ricci equations
\begin{eqnarray}
d\theta_2^3 + \theta_1^3 \wedge \theta_2^1 + \theta_4^3 \wedge \theta_2^4 + \theta_5^3 \wedge \theta_2^5  &  =  & 0 \nonumber \\
d\theta_1^4 + \theta_2^4 \wedge \theta_1^2 + \theta_3^4 \wedge \theta_1^3 + \theta_5^4 \wedge \theta_1^5  &  =  & 0 \nonumber \\
d\theta_3^5 + \theta_1^5 \wedge \theta_3^1 + \theta_2^5 \wedge \theta_3^2 + \theta_4^5 \wedge \theta_3^4  &  =  & 0 \nonumber
\end{eqnarray}
We get
\begin{eqnarray}\label{eq:Codazzi2}
&& (a_1+b_2) + b \cot\alpha \alpha_2 + a (\cot\alpha \alpha_1 + 6 \tan\beta
\cos\alpha)\nonumber \\
&& - 2 \sec\beta \cos\alpha \alpha_2=0
\end{eqnarray}
Using the connection form (\ref{eq:normalconexaobeta}) in the Codazzi-Ricci equations 
\begin{eqnarray}
d\theta_2^4 + \theta_1^4 \wedge \theta_2^1 + \theta_3^4 \wedge \theta_2^3 + \theta_5^4 \wedge \theta_2^5  &  =  & 0 \nonumber \\
d\theta_4^5 + \theta_1^5 \wedge \theta_4^1 + \theta_2^5 \wedge \theta_4^2 + \theta_3^5 \wedge \theta_4^3  &  =  & 0 \nonumber \\
d\theta_3^4 + \theta_1^4 \wedge \theta_3^1 + \theta_2^4 \wedge \theta_3^2 + \theta_5^4 \wedge \theta_3^5  &  =  & 0 \nonumber
\end{eqnarray}
We have
\begin{eqnarray}\label{eq:Codazzi3}
&& (a_2-b_1) - (a^2+ b^2) \cot\alpha \sin\beta \cot^2 \beta + a \cot\alpha \alpha_2  \\
&& + b( - cot\alpha \alpha_1 + 2 \cos\alpha (\cot\beta -3 \tan\beta)) + 2 \cos\alpha \sin\beta (\cot\beta -\tan\beta)\alpha_1 \nonumber \\
&& + \sin\alpha \cos\alpha \sin\beta (5 - \cot^2\beta) + \sin \beta \Delta \alpha  =  0\nonumber
\end{eqnarray}
Codazzi-Ricci equations
\begin{eqnarray}
d\theta_1^2 + \theta_3^2 \wedge \theta_1^3 + \theta_4^2 \wedge \theta_1^4 + \theta_5^2 \wedge \theta_1^5  &  =  & 0 \nonumber \\
d\theta_1^5 + \theta_2^5 \wedge \theta_1^2 + \theta_3^5 \wedge \theta_1^3 + \theta_4^5 \wedge \theta_1^4  &  =  & 0 \nonumber
\end{eqnarray}
give the following equation
\begin{eqnarray}\label{eq:Codazzi4}
&& (a^2 + b^2) (1+\csc^2\beta) + 2b\csc\beta(\alpha_1 - \cot\beta \sin\alpha)
-2a\csc\beta\alpha_2  \nonumber\\
&&  + |\nabla \alpha|^2 + 2\sin\alpha (\tan \beta - \cot\beta) \alpha_1 - 4
\tan^2 \beta \cos^2 \alpha \nonumber \\
&& -\sin^2 \alpha (1-\cot^2\beta) = 0
\end{eqnarray}
The following Codazzi equation is automatically verified
\begin{eqnarray}
d\theta_2^5 + \theta_1^5 \wedge \theta_2^1 + \theta_3^5 \wedge \theta_2^3 + \theta_4^5 \wedge \theta_2^4  &  =  & 0 \nonumber
\end{eqnarray}
\section{Aplications}
\subsection{Minimal Surfaces in $S^5$  with Constant Contact Angle $\beta$ and parallel normal vector field} \mbox{}

In this section, we will give Gauss-Codazzi-Ricci equations for a minimal surface in $S^5$ with constant contact angle and null principal curvatures $a,b$.\\
Syppose that $a,b$ are nulls and the Contact angle $\beta$ is constant, then using the  Codazzi equation (\ref{eq:Codazzi1}), we have \begin{eqnarray}\label{eq:Codred1}
\cos\alpha (2(\cot\beta - \tan\beta)\alpha_1  - \sin\alpha (\cot^2 \beta -3)) - \cot\alpha |\nabla \alpha|^2 =  0
\end{eqnarray}
On the other hand, Codazzi equation (\ref{eq:Codazzi3}) with $a,b$ nulls and constant contact angle implies
\begin{eqnarray}\label{eq:Codred2}
2 \cos\alpha (\cot\beta -\tan\beta)\alpha_1 +  \sin\alpha \cos\alpha  (5 - \cot^2\beta) +  \Delta \alpha  =  0
\end{eqnarray}
Using equations (\ref{eq:Codred1}) and (\ref{eq:Codred2}), we obtain a new Laplacian equation of $\alpha$
\begin{eqnarray}
\Delta \alpha & = & -\sin(2\alpha) -\cot\alpha  |\nabla \alpha|^2
\end{eqnarray}
Now using the Hopf`s Lemma, we get the Theorem 1.
\begin{observacao}
The Theorem 1 implies a general classification in \cite{RMV} that gives a family of minimal flat tori in $S^5$ with constant Contact angle  and constant Holomorphic angle 
\end{observacao}

\end{document}